\documentclass[12pt, reqno]{amsart}
\usepackage{amsmath, amsthm, amscd, amsfonts, amssymb, graphicx, color, enumerate}
\usepackage[bookmarksnumbered, colorlinks, plainpages]{hyperref}
\usepackage[all,cmtip]{xy}
\usepackage{mathrsfs}
\hypersetup{colorlinks=true,linkcolor=red, anchorcolor=green, citecolor=cyan, urlcolor=red, filecolor=magenta, pdftoolbar=true}

\allowdisplaybreaks[4]

\newtheorem{theorem}{Theorem}[section]
\newtheorem{lemma}[theorem]{Lemma}

\newtheorem{corollary}[theorem]{Corollary}
\theoremstyle{definition}
\newtheorem{definition}[theorem]{Definition}
\newtheorem{example}[theorem]{Example}

\newtheorem{conjecture}[theorem]{Conjecture}

\newtheorem{problem}[theorem]{Problem}
\newtheorem{remark}[theorem]{Remark}
\numberwithin{equation}{section}

\begin{document}
\setcounter{page}{1}

\title[Compactness of composition operators on Hardy spaces]{Compactness of composition operators on Hardy spaces of homogeneous rooted trees}

\author[L. Liu]{Lei Liu$^{*}$}

\address{School of Mathematics and Statistics, Xidian University, Xi'an 710126, P.R. China.}

\email{\textcolor[rgb]{0.00,0.00,0.84}{leiliu@mail.xidian.edu.cn}}

\author[W. Ma]{Wanxin Ma}

\address{School of Mathematics and Statistics, Xidian University, Xi'an 710126, P.R. China.}

\email{\textcolor[rgb]{0.00,0.00,0.84}{mawanxin26@126.com}}
%\dedicatory{This paper is dedicated to Professor ABCD}

\makeatletter
\@namedef{subjclassname@2020}{\textup{2020} Mathematics Subject Classification}
\makeatother
\subjclass[2020]{Primary: 47B33; Secondary: 05C05.}

\keywords{Composition operator; Compact operator; Linear combination; Discrete Hardy space; Homogeneous rooted tree}

\date{Received: xxxxxx; Revised: yyyyyy; Accepted: zzzzzz.
\newline \indent $^{*}$ Corresponding author}

\begin{abstract}
In this paper, we present a complete characterization of the compactness of arbitrary finite linear combinations of composition operators acting on Hardy spaces of homogeneous rooted trees. In particular, we obtain a characterization of compactness for linear combinations satisfying the CNC. As an application, we give an affirmative answer to a conjecture raised by Choe et al. [Trans. Amer. Math. Soc. Ser. B 9(2022), 733--756] concerning such Hardy spaces. Furthermore, we establish the phenomenon of (double) difference cancellation for linear combinations of composition operators by constructing explicit examples.
\end{abstract} \maketitle

\section{Introduction}

Let $X$ be a Banach space consisting of complex-valued functions on a set $\Omega$. For a self-map $\varphi$ of $\Omega$, the composition operator $C_{\varphi}$ is defined as
\[
C_{\varphi}(f)=f\circ\varphi \ \ \text{for}\ \ f\in X.
\]
Composition operators have been extensively investigated on various Banach spaces of holomorphic functions over the past half century. The main subject in the study of composition operators is to describe operator theoretic properties of $C_{\varphi}$ in terms of function theoretic properties of $\varphi$ (see \cite{ba,cm,do,sh,zh}).

One central problem in the study of composition operators lies in characterizing the compactness of these operators as well as the compactness of their differences. Shapiro \cite{sh2} established a full characterization for the compactness of composition operators acting on Hardy spaces over the unit disk. MacCluer and Shapiro \cite{ms} derived necessary and sufficient conditions for compact composition operators on weighted Bergman spaces in terms of angular derivatives. In \cite{wzz}, Wulan et al. provides a criterion for compactness of composition operator $C_\varphi$ in terms of the vanishing behaviour of the norms of the powers of the symbol function $\varphi$ on BMOA and Bloch spaces. Moorhouse \cite{mo} characterized compact differences of composition operators on weighted Bergman spaces by the Julia-Carath\'eodory type condition and found that two noncompact composition operators can induce compact differences of composition operators, which is called Moorhouse cancellation property. In \cite{ckp,ckp2} this characterization was extended to higher dimensional polydisks and balls.

In recent years, numerous papers have investigated the compactness of composition operators defined on discrete structures including trees. Allen, Colonna, and Easley \cite{ace} obtained a characterization for the compactness of composition operators acting on the Lipschitz space of a tree. In \cite{lk} Liu and Kong characterized compact differences of composition operators on the Lipschitz space of a tree and further verified that the Moorhouse cancellation property holds on this Lipschitz space.
Muthukumar and Ponnusamy \cite{mp} studied the compactness of composition operators acting on Hardy spaces of homogeneous rooted trees and obtained a sufficient condition for such operators to be compact.
Motivated by the above results, we pose the following problem:
\begin{problem}\label{pro}
\textit{Find necessary and sufficient conditions that characterize the compactness of composition operators acting on Hardy spaces of the homogeneous rooted trees.}
\end{problem}
The main purpose of the present paper is to establish a complete characterization for compact linear combinations of a finite number of composition operators acting on Hardy spaces of homogeneous rooted trees. As an immediate corollary, we resolve the aforementioned problem.

The paper is organized as follows. We refer to Section 2 for some preliminaries on trees and the definitions of Hardy spaces of homogeneous rooted trees. In Section 3, we characterize the compactness of arbitrary finite linear combinations of composition operators acting on these Hardy spaces. As a direct consequence of our main result, we settle Problem \ref{pro}. Furthermore, we obtain a characterization of the compactness of linear combinations satisfying the CNC and confirm a conjecture due to Choe et al. In Section 4, we construct two concrete examples to illustrate the phenomenon of (double) difference cancellation on Hardy spaces over homogeneous rooted trees.

\section{preliminaries}

A tree is a locally finite connected and simply connected graph $T$, which, as a set, we identify with the collection of its vertices. Two vertices $u$ and $v$ are said to be neighbors if there is an edge connecting them, and we adopt the notation $u \sim v$. The degree of a vertex is the number of its neighbors. A tree is said to be $k$-homogeneous if every vertex has degree exactly $k$. A path between two vertices $u$ and $w$ is a finite sequence of distinct vertices $[u, u_1, u_2, \ldots, u_{k}, w]$ such that $u \sim u_1 \sim u_2 \sim \cdots \sim u_{k} \sim w$. Since $T$ is a tree, the path between two vertices is unique.
	
Given a tree $T$ rooted at $o$ and a vertex $u\in T$, a vertex $v \neq u$ is called a descendant of $u$ if $u$ belongs to the unique path from $o$ to $v$. In this case, $u$ is called an ancestor of $v$. We write $v^-$ for the parent of $v$, namely the unique neighbor of $v$ which is an ancestor of $v$. The vertex $v$ is called a child of $v^-$. The length of a path between two vertices $u$ and $w$ is the number of edges connecting them. The distance $d(u, w)$ between vertices $u$ and $w$ is the length of the path connecting $u$ to $w$. In particular, the length of a path joining the root $o$ to a vertex $w$ is called the length of $w$ and is denoted by $|w|$.	
A function on a tree is a complex-valued function on the set of its vertices.

Throughout this paper, $T$ always stands for a $(q+1)$-homogeneous rooted tree with root vertex $o$, and we adopt the convention $\mathbb{N}_{0}=\mathbb{N} \cup \{0\}$. For each $n\in\mathbb{N}$, we define
\[
M_{p}(n, f):= \begin{cases}\left(\frac{1}{(q+1) q^{n-1}} \sum_{|v|=n}|f(v)|^{p}\right)^{\frac{1}{p}} & \text{if } p \in(0, \infty), \\ \max _{|v|=n}|f(v)| & \text { if } p=\infty \end{cases}
\]
and
\[
\| f\| _{p}:=\sup _{n \in \mathbb{N}_{0}} M_{p}(n, f),
\]
with the convention $M_{p}(0, f):=|f(o)|$.
The discrete Hardy space associated with the tree $T$, denoted by $\mathbb{T}_{p}$, is defined as
\[
\mathbb{T}_{p}:=\left\{f: T \to \mathbb{C} \mid\| f\| _{p}<\infty\right\}
\]
for all $p \in(0, \infty]$. By \cite[Theorem A]{mp}, $\mathbb{T}_{p}$ forms a Banach space for $1\leq p\leq\infty$. For the sake of simplicity, we omit the parameter $q$ from the notation for this space.
	
We fix the following notation for the remainder of the paper. For each $n \in \mathbb{N}_{0}$, write $D_{n}$ for the collection of vertices $v \in T$ satisfying $|v|=n$, and let $c_{n}$ denote the cardinality of $D_{n}$. Accordingly,
\[
c_{n}= \begin{cases}(q+1) q^{n-1} & \text { if } n \in \mathbb{N}, \\ 1 & \text { if } n=0. \end{cases}
\]
For $n\in\mathbb{N}_{0}$ and $w\in T$, let $N_{\varphi}(n,w)$ stand for the number of preimages of $w$ under $\varphi$ lying in the set $D_{n} = \{v\in T:|v|=n\}$. Equivalently,
\[
N_{\varphi}(n,w) = \#\big(\varphi^{-1}(w)\cap D_{n}\big).
\]
Furthermore, for all $m,n\in\mathbb{N}_{0}$, define $N_{m,n}$ as the maximum of $N_{\varphi}(n,w)$ taken over all vertices $w$ with $|w|=m$. It is straightforward to verify that
\[
\sum_{m=0}^{\infty} N_{m,n} \leq c_{n}
\]
holds for every nonnegative integer $n$.

\section{Results and proofs}

In this section, we establish a necessary and sufficient condition for the compactness of arbitrary finite linear combinations of composition operators acting on $\mathbb T_p$. We first present a useful compactness criterion for such linear combinations on $\mathbb T_p$. A proof for the case of a single composition operator is available in \cite[Theorem G]{mp}. For the sake of completeness, we supply a detailed proof for the linear combination case below.

\begin{lemma}\label{lem}	
Let $T$ be a $(q+1)$-homogeneous rooted tree with $q\geq2$, and let $1\leq p<\infty$. For each $i=1, 2, \dots, m$, suppose that $\varphi_i$ are self-maps of $T$ and $C_{\varphi_i}$ denote the composition operators on $\mathbb T_p$. Define
\[
Q:=\sum_{i=1}^{m}a_iC_{\varphi_i},
\qquad a_i\in\mathbb C\setminus\{0\}.
\]
Then the following statements are equivalent:
\begin{enumerate}[(i)]
\item $Q$ is compact on $\mathbb T_p$;
\item $\|Qf_j\|\to0$ for every bounded sequence $\{f_j\}_{j=1}^{\infty}\subset\mathbb T_p$ that converges pointwise to $0$.
\end{enumerate}
\end{lemma}

\begin{proof}
$``(i)\Rightarrow(ii)"$	 Assume that $\{f_j\}_{j=1}^{\infty}\subset\mathbb T_p$	is a bounded sequence converging pointwise to $0$. If $\|Qf_j\|\nrightarrow0$, then there exist a subsequence $\{f_{j_k}\}_{k=1}^{\infty}$ and a constant $\varepsilon>0$ satisfying
\[
\|Qf_{j_k}\|\geq\varepsilon
\]
for all $k$. Since $Q$ is compact, there exists a further subsequence, still denoted by $\{f_{j_k}\}_{k=1}^{\infty}$, such that $\{Qf_{j_k}\}_{k=1}^{\infty}$ converges in norm
to some $g\in\mathbb T_p$. On the one hand, since $\{f_{j_k}\}_{k=1}^{\infty}$ converges pointwise to $0$, we have
\[
Qf_{j_k}(v)
=\sum_{i=1}^{m}a_if_{j_k}(\varphi_i(v))
\rightarrow0,
\]
for each fixed vertex $v$ in $T$ as $k\to\infty$.
On the other hand, norm convergence in
$\mathbb T_p$ trivially implies pointwise convergence. Consequently, $g(v)=0$ for all $v$ in $T$, so that $g\equiv0$. It follows that
\[
\|Qf_{j_k}\|\rightarrow0,
\]
as $k\to\infty$, which contradicts $\|Qf_{j_k}\|\geq\varepsilon$.
		
$``(ii)\Rightarrow(i)"$	
Let $\{g_j\}_{j=1}^{\infty}$ be a sequence in the unit ball of $\mathbb T_p$. For each fixed vertex $v\in T$, we have
\[
|g_j(v)|\leq c_{|v|}^{1/p}M_{p}(|v|, g_j)\leq c_{|v|}^{1/p},
\]
which implies $\{g_j(v)\}_{j=1}^{\infty}$ is bounded in $\mathbb{C}$. Applying a diagonalization argument, we may extract a subsequence $\{g_{j_k}\}_{k=1}^{\infty}$ that converges pointwise to some function $g: T\to\mathbb C$.
		
Moreover, for each $n\in\mathbb N_0$, we have
\[
M_p^p(n,g)
=\lim_{k\to\infty}\frac1{c_n}
\sum_{|v|=n}|g_{j_k}(v)|^p
\leq1.
\]
Therefore $g\in\mathbb T_p$ and $\|g\|\leq1$. Set $f_k=g_{j_k}-g$.
Then the sequence $\{f_k\}_{k=1}^{\infty}$ is bounded in $\mathbb T_p$ and converges pointwise to $0$. By condition (ii), we obtain
\[
\|Qf_k\|\rightarrow 0,
\]
as $k\to\infty$, which is equivalent to
\[
\|Qg_{j_k}-Qg\|\rightarrow 0.
\]
Thus, every sequence contained in the image of the unit ball of $\mathbb T_p$ under $Q$ admits a subsequence convergent with respect to the norm. Therefore, $Q$ is a compact operator.
\end{proof}

We now introduce the following definition, which will play a key role throughout this paper.

\begin{definition}
Given
\[
Q=\sum_{i=1}^{r}a_iC_{\varphi_i},
\qquad a_i\in\mathbb C\setminus\{0\},
\]
where $\varphi_i:T\to T$, $i=1,\ldots,r$, are self-maps. For any vertices $v,w\in T$, we define
\[
A_Q(v,w):=(Q\chi_{\{w\}})(v)
=\sum_{i=1}^r a_i\chi_{\{w\}}(\varphi_i(v)).
\]
For any $n\in\mathbb N_0$ and $w\in T$, we further define
\[
N_Q(n,w):=\sum_{|v|=n}|A_Q(v,w)|^p.
\]
\end{definition}

Our main result is the next theorem.

\begin{theorem}\label{thm}
Let $T$ be a $(q+1)$-homogeneous tree with $q\ge2$, and let $1\le p<\infty$. For each $i=1, \ldots, r$, let $\varphi_i: T\to T$ be a self-map such that the composition operator $C_{\varphi_i}$ is bounded on $\mathbb T_p$. Then the operator $Q=\sum_{i=1}^{r} a_iC_{\varphi_i}$, $0\ne a_i\in\mathbb C$, is compact on $\mathbb T_p$ if and only if
\[
\lim_{R\to\infty}\{\sup_{n\in\mathbb N_0}
\frac1{c_n}\sum_{l>R}c_l\max_{|w|=l}N_Q(n,w)\}=0.
\]
\end{theorem}
	
\begin{proof}
$``\Leftarrow"$ Let $\{f_j\}_{j=1}^{\infty}$ be a bounded sequence in $\mathbb T_p$ converging pointwise to $0$. Without loss of generality, assume $\|f_j\|\le1$ for all $j\in\mathbb{N}$. It suffices to prove that $\|Qf_j\|\to0$ as $j\to\infty$. Since $\|f_j\|\le1$, we have
\begin{equation}
\sum_{|w|=l}|f_j(w)|^p\le c_l
\end{equation}
for all $j,l\in\mathbb N$. Set $A=\sum_{i=1}^{r} |a_i|$. For any fixed vertex $v\in T$, there exist at most $r$ vertices $w\in T$ satisfying $A_Q(v,w)\ne0$, and $|A_Q(v,w)|\le A$. Hence, for any $n\in\mathbb N_0$, we have
\begin{equation}
\sum_{w\in T}N_Q(n,w)
=\sum_{|v|=n}\sum_{w\in T}|A_Q(v,w)|^p
\le rA^pc_n.
\end{equation}
For any $\varepsilon>0$, we may choose $R\in\mathbb N_0$ such that
\begin{equation}
\sup_{n\in\mathbb N_0}
\frac1{c_n}\sum_{l>R}c_l\max_{|w|=l}N_Q(n,w)
<\frac{\varepsilon^p}{2r^{p-1}}.
\end{equation}
Set $S=\bigcup_{l=0}^{R}D_l$. Since $S$ is finite and $f_j\to0$ pointwise, there exists $J\in\mathbb N$ satisfying
\begin{equation}
\sup_{w\in S}|f_j(w)|^p<\frac{\varepsilon^p}{2r^pA^p}
\end{equation}
whenever $j\ge J$.
Moreover, for any $j\ge J$ and $n\in\mathbb N_0$, we have
\begin{align*}
M_p^p(n,Qf_j)
&=\frac1{c_n}\sum_{|v|=n}|Qf_j(v)|^p  \\
&=\frac1{c_n}\sum_{|v|=n}
\left|\sum_{i=1}^{r} a_i
\sum_{w\in T}\chi_{\{w\}}(\varphi_i(v))f_j(w)\right|^p \\
&=\frac1{c_n}\sum_{|v|=n}
\left|\sum_{w\in T}A_Q(v,w)f_j(w)\right|^p \\
&\le \frac{r^{p-1}}{c_n}\sum_{|v|=n}
\sum_{w\in T}|A_Q(v,w)|^p|f_j(w)|^p \\
&=\frac{r^{p-1}}{c_n}\sum_{w\in T}
\left(\sum_{|v|=n}|A_Q(v,w)|^p\right)|f_j(w)|^p \\
&=\frac{r^{p-1}}{c_n}\sum_{w\in T}N_Q(n,w)|f_j(w)|^p .
\end{align*}
Define
\[
I:=\frac{r^{p-1}}{c_n}\sum_{w\in S}N_Q(n,w)|f_j(w)|^p
\]
and
\[
II:=\frac{r^{p-1}}{c_n}\sum_{l>R}\sum_{|w|=l}N_Q(n,w)|f_j(w)|^p.
\]
It follows that		
\begin{align*}
M_p^p(n,Qf_j)\leq I+II.
\end{align*}
Using (3.2) and (3.4), we deduce that
\begin{align*}
I&\le \frac{r^{p-1}}{c_n}\left(\sum_{w\in T}N_Q(n,w)\right)\sup_{w\in S}|f_j(w)|^p  \\
&\le r^pA^p\sup_{w\in S}|f_j(w)|^p
<\frac{\varepsilon^p}{2}.
\end{align*}
Moreover, since $\|f_j\|\le1$, it follows from (3.1) that
\[
\sum_{|w|=l}N_Q(n,w)|f_j(w)|^p
\le c_l\max_{|w|=l}N_Q(n,w).
\]
Combining the above inequality with (3.3), we therefore obtain that
\[
II\le r^{p-1}\left(\frac1{c_n}\sum_{\ell>R}c_\ell\max_{|w|=\ell}N_Q(n,w)\right)<\frac{\varepsilon^p}{2}.
\]
From all the preceding estimates, we conclude that $M_p^p(n,Qf_j)<\varepsilon^p$ holds for all $n\in\mathbb N_0$ and $j\ge J$. Hence $\|Qf_j\|<\varepsilon$ for all $j\ge J$, which implies $\|Qf_j\|\to0$ as $j\to\infty$. By Lemma \ref{lem}, the operator $Q$ is compact on $\mathbb T_p$.
		
$``\Rightarrow"$ We proceed by contradiction. Assume that there exist $\delta>0$, a sequence $R_j\to\infty$, and $n_j\in\mathbb N_0$ such that
\[
\frac1{c_{n_j}}\sum_{l>R_j}c_l\max_{|w|=l}N_Q(n_j,w)>\delta
\]
for all $j$. For each $l>R_j$, choose $w_{j,l}\in D_l$ such that
\[
N_Q(n_j,w_{j,l})=\max_{|w|=l}N_Q(n_j,w).
\]
Set
\[
\Lambda_j=\{l>R_j: \exists ~v\in D_{n_j} ~s.t. \ A_Q(v,w_{j,l})\ne0\},
\qquad L_j=\#\Lambda_j.
\]
It follows that $\Lambda_j$ is a finite set. Indeed, $D_{n_j}$ is finite, and $A_Q(v,w_{j,l})\ne0$ only if $w_{j,l}\in\{\varphi_1(v), \ldots, \varphi_r(v)\}$. In particular, for each fixed $v\in D_{n_j}$, there are at most $r$ indices $l>R_j$ such that $A_Q(v,w_{j,l})\ne0$.

We proceed by establishing an inequality: for any $z_1, \ldots, z_k\in\mathbb C$, there exists a constant $\eta_{p,k}>0$ such that
\[
\frac1{2^k}\sum_{\delta_1,\ldots,\delta_k=\pm1}
|\delta_1z_1+\cdots+\delta_kz_k|^p
\ge \eta_{p,k}\sum_{s=1}^k|z_s|^p.
\]
Without loss of generality, we may assume $\sum_{s=1}^k|z_s|^p=1$. On the compact set $\{(z_1, \ldots, z_k)\in\mathbb{C}^k: \sum_{s=1}^k|z_s|^p=1$\}, the function
\[
\frac1{2^k}\sum_{\delta_1,\ldots,\delta_k=\pm1}
|\delta_1z_1+\cdots+\delta_kz_k|^p
\]
is continuous and hence attains its minimum. This minimum is strictly positive. Indeed, suppose it were zero. Then there exists some $(z_{1}^{(0)},\ldots,z_{k}^{(0)})$ satisfying $\sum_{s=1}^k|z_{s}^{(0)}|^p=1$, such that
\[
\frac1{2^k}\sum_{\delta_1,\ldots,\delta_k=\pm1}
|\delta_1z_{1}^{(0)}+\cdots+\delta_kz_{k}^{(0)}|^p=0,
\]
which implies $\delta_1z_{1}^{(0)}+\cdots+\delta_kz_{k}^{(0)}=0$ for all $\delta_1,\ldots,\delta_k=\pm1$.
Fixing $s$ and comparing two sign choices which differ only in the $s$ sign, we see $2z_{s}^{(0)}=0$. Hence $z_{s}^{(0)}=0$ for all $s\in\{1,\ldots,k\}$. This is a contradiction. Therefore there exists $\eta_{p,k}>0$ such that
\[
\frac1{2^k}\sum_{\delta_1,\ldots,\delta_k=\pm1}
|\delta_1z_1+\cdots+\delta_kz_k|^p
\ge
\eta_{p,k}\sum_{s=1}^k|z_s|^p.
\]
Let $\eta_{p,r}:=\min_{1\leq k\leq r}\{\eta_{p,k}\}>0$. This is a constant depending only on $p$ and $r$. It follows that
\[
\frac1{2^k}\sum_{\delta_1,\ldots,\delta_k=\pm1}
|\delta_1z_1+\cdots+\delta_kz_k|^p
\ge \eta_{p,r}\sum_{s=1}^k|z_s|^p.
\]
		
Next, notice that for each fixed $v\in D_{n_j}$, there are at most $r$ indices $l\in\Lambda_j$ satisfying $A_Q(v,w_{j,l})\ne0$. Applying the preceding estimate to the nonzero elements of $\{c_l^{1/p}A_Q(v,w_{j,l}): l\in\Lambda_j\}$, we deduce that
\[
\frac1{2^{L_j}}\sum_{\varepsilon_{j,l}=\pm1,\ l\in\Lambda_j}
\left|\sum_{l\in\Lambda_j}\varepsilon_{j,l}c_l^{1/p}A_Q(v,w_{j,l})\right|^p
\ge\eta_{p,r}\sum_{l\in\Lambda_j}c_l|A_Q(v,w_{j,l})|^p.
\]
Taking the sum over all $v\in D_{n_j}$, we arrive at
\[
\begin{aligned}
&\frac1{2^{L_j}}\sum_{\varepsilon_{j,l}=\pm1,l\in\Lambda_j}\sum_{|v|=n_j}
\left|\sum_{l\in\Lambda_j}\varepsilon_{j,l}c_l^{1/p}A_Q(v,w_{j,l})\right|^p \\
&\quad\ge\eta_{p,r}\sum_{|v|=n_j}\sum_{l\in\Lambda_j}c_l|A_Q(v,w_{j,l})|^p.
\end{aligned}
\]
Therefore, there exists a choice of signs $\varepsilon_{j,l}\in\{-1,1\}$ for all $l\in\Lambda_j$ such that
\[
\sum_{|v|=n_j}
\left|\sum_{l\in\Lambda_j}\varepsilon_{j,l}c_l^{1/p}A_Q(v,w_{j,l})\right|^p
\ge\eta_{p,r}\sum_{|v|=n_j}\sum_{l\in\Lambda_j}c_l|A_Q(v,w_{j,l})|^p.
\]
Set $\varepsilon_{j,l}=1$ for all $l\notin\Lambda_j$. As $A_Q(v,w_{j,l})=0$ whenever $v\in D_{n_j}$ and $l\notin\Lambda_j$, the preceding inequality simplifies to
\begin{equation}
\sum_{|v|=n_j}
\left|\sum_{l>R_j}\varepsilon_{j,l}c_l^{1/p}A_Q(v,w_{j,l})\right|^p
\ge\eta_{p,r}\sum_{l>R_j}c_l N_Q(n_j,w_{j,l}).
\end{equation}
Now define
\[
f_j:=\sum_{l>R_j}\varepsilon_{j,l}c_l^{1/p}\chi_{\{w_{j,l}\}}.
\]
Then $\|f_j\|=1$ for all $j\in\mathbb{N}$, and $f_j\to0$ pointwise. For each $v\in D_{n_j}$, by the definition of $A_Q$, we have
\[
Qf_j(v)=\sum_{i=1}^{r} a_if_j(\varphi_i(v))
=\sum_{l>R_j}\varepsilon_{j,l}c_l^{1/p}A_Q(v,w_{j,l}).
\]
Consequently, it follows from (3.5) that
\[
\begin{aligned}
M_p^p(n_j,Qf_j)
&=\frac1{c_{n_j}}\sum_{|v|=n_j}|Qf_j(v)|^p\\
&=\frac1{c_{n_j}}\sum_{|v|=n_j}
\left|\sum_{l>R_j}\varepsilon_{j,l}c_l^{1/p}A_Q(v,w_{j,l})\right|^p\\
&\ge \frac{\eta_{p,r}}{c_{n_j}}\sum_{l>R_j}c_l N_Q(n_j,w_{j,l})\\
&=\frac{\eta_{p,r}}{c_{n_j}}\sum_{l>R_j}c_l\max_{|w|=l}N_Q(n_j,w)\\
&>\eta_{p,r}\delta,
\end{aligned}
\]
which implies
\[
\|Qf_j\|>(\eta_{p,r}\delta)^{1/p}>0
\]
for all $j\in\mathbb{N}$. However, by Lemma \ref{lem}, the compactness of $Q$ ensures that $\|Qf_j\|\to0$ as $j\to\infty$, which yields a contradiction. This completes the proof.
\end{proof}

The following corollary gives a necessary and sufficient characterization of compact composition operators on $\mathbb{T}_p$, thereby settling Problem \ref{pro}.

\begin{corollary}\label{cor3.4}
Let $T$ be a $(q+1)$-homogeneous rooted tree with $q\geq 2$ and $1\leq p<\infty$. Let $\varphi$ be a self-map of $T$ such that $C_{\varphi}$ is bounded. Then composition operator $C_{\varphi}:\mathbb{T}_p\to\mathbb{T}_p$ is compact if and only if
\[
\lim_{R\to\infty}\{\sup_{n\in\mathbb{N}_0}
\frac{1}{c_n}\sum_{m>R}c_mN_{m,n}\}=0.
\]
\end{corollary}
	
\begin{proof}
By Theorem \ref{thm}, it suffices to verify that $N_{m,n}=\max_{|w|=m}N_{C_{\varphi}}(n,w)$
for all $m, n\in\mathbb{N}_0$. Indeed, for any $v, w\in T$, by the definition of
$A_{C_{\varphi}}$, we have
\[
A_{C_{\varphi}}(v,w)
=(C_{\varphi}\chi_{\{w\}})(v)
=\chi_{\{w\}}(\varphi(v)).
\]
Hence
\begin{align*}
N_{C_{\varphi}}(n,w)
&=\sum_{|v|=n}|A_{C_{\varphi}}(v,w)|^p\\
&=\sum_{|v|=n}
\left|\chi_{\{w\}}(\varphi(v))\right|^p\\
&=\#\{v\in D_n:\varphi(v)=w\}\\
&=N_{\varphi}(n,w).
\end{align*}
From this we deduce that for each $m,n\in\mathbb{N}_0$,
\[
\max_{|w|=m}N_{C_{\varphi}}(n,w)
=\max_{|w|=m}N_{\varphi}(n,w)
=N_{m,n}
\]
as desired.
\end{proof}

Now we consider the class of linear combinations satisfying the coefficient non-cancellation condition. Let $\mathscr{P}_r$ denote the group of all permutations on the index set $\{1,\ldots,r\}$.
Identifying $\sigma\in\mathscr{P}_r$ with the ordered $r$-tuple $(\sigma_1,...,\sigma_r)$ where $\sigma_i:=\sigma(i)$, we may write $\sigma=(\sigma_1,\ldots,\sigma_r)$. For any $\sigma\in\mathscr{P}_r$, define
\[
s_i^\sigma:=\sum_{k=1}^{i}a_{\sigma_k},
\qquad 1\leq i\leq r.
\]
Given a linear combination $Q=\sum_{i=1}^{r}a_iC_{\varphi_i}$, we say that $Q$ satisfies the coefficient non-cancellation condition (CNC) provided
\[
s_j^\sigma\ne0
\ \ \text{for all} \ \
1\leq j\leq r-1 \ \ \text{and} \ \ \sigma\in \mathscr{P}_r.
\tag{CNC}
\]
The following characterization holds for such linear combinations.

\begin{theorem}\label{thm3.5}
Let $T$ be a $(q+1)$-homogeneous tree with $q\ge2$, and let $1\le p<\infty$. For each $i=1, \ldots, r$, let $\varphi_i: T\to T$ be a self-map such that the composition operator $C_{\varphi_i}$ is bounded on $\mathbb T_p$. Assume that the operator $Q=\sum_{i=1}^{r} a_iC_{\varphi_i}$, $0\ne a_i\in\mathbb C$ satisfies the CNC. Then $Q$ is compact on $\mathbb T_p$ if and only if one of the following holds.
\begin{enumerate}[(i)]
\item $\sum_{i=1}^{r}a_i\ne0$, and $C_{\varphi_i}$ is compact on $\mathbb T_p$ for all $i=1,\ldots,r$;
\item $\sum_{i=1}^{r}a_i=0$, and $C_{\varphi_i}-C_{\varphi_j}$ is compact on $\mathbb T_p$ for all $i,j\in\{1,\ldots,r\}$.
\end{enumerate}
\end{theorem}

\begin{proof}
The sufficiency follows immediately for case (i).

For case (ii), since $\sum_{i=1}^{r}a_i=0$, we have
\[
\sum_{i=1}^{r}a_iC_{\varphi_i}=\sum_{i=1}^{r-1}
a_i(C_{\varphi_i}-C_{\varphi_r}).
\]
By our assumption, each $C_{\varphi_i}-C_{\varphi_r}$ is compact on $\mathbb T_p$. Therefore, $\sum_{i=1}^{r}a_iC_{\varphi_i}$ is compact.
    	
We now turn to the proof of necessity. Suppose
\[
\sum_{i=1}^{r}a_i\ne0.
\]
For any $i=1,\ldots,r$,	we will show that $C_{\varphi_i}$ is compact. Indeed, set
\[
\xi=\min_{\substack{\sigma\in \mathscr{P}_r\\1\leq j\leq r}}|s_j^\sigma|^p.
\]
Since the coefficients satisfy the CNC and $\sum_{i=1}^{r}a_i\ne0$, it follows that $\xi>0$.
For any $v, w\in T$, if $\varphi_k(v)=w$, then put $I=\{k\in\{1,\ldots,r\}:\varphi_k(v)=w\}$.
We can select $\sigma=(\sigma_1,\ldots,\sigma_r)\in \mathscr{P}_r$ so that $\sigma_1,\ldots,\sigma_{|I|}$ are exactly the elements of $I$. It follows that
\begin{align*}
A_Q(v,w)&=\sum_{k=1}^{r}a_k\chi_{\{w\}}(\varphi_k(v))\\
&=\sum_{l=1}^{|I|}a_{\sigma_l}=s_{|I|}^{\sigma}.
\end{align*}
Thus,
\[
|A_Q(v,w)|^p\geq \xi.
\]
Without lose of generality, we may assume that $|v|=n$ and $|w|=m$. Then
\begin{align*}
N_Q(n,w)
&=\sum_{|v|=n}|A_Q(v,w)|^p\\
&\geq\xi\,\#\{v\in D_n:\varphi_i(v)=w\}\\
&=\xi N_{\varphi_i}(n,w),
\end{align*}
and hence
\[
\max_{|w|=m}N_{\varphi_i}(n,w)
\leq\frac1{\xi}\max_{|w|=m}N_Q(n,w).
\]
Since $Q$ is compact, it follows from Theorem \ref{thm} that
\begin{align*}
0&\leq\sup_{n\in\mathbb N_0}\frac1{c_n}\sum_{m>R}c_m\max_{|w|=m}N_{\varphi_i}(n,w)\\
&\leq\frac1{\xi}\sup_{n\in\mathbb N_0}\frac1{c_n}\sum_{m>R}c_m\max_{|w|=m}N_Q(n,w)\to0
\end{align*}
as $R\to\infty$. Applying Corollary \ref{cor3.4}, we conclude that $C_{\varphi_i}$ is compact on $\mathbb T_p$, as desired.
    	
Next, suppose
\[
\sum_{i=1}^{r}a_i=0.
\]
For any distinct $i,j\in\{1,\ldots,r\}$, set
$D=C_{\varphi_i}-C_{\varphi_j}$.
We shall show that $D$ is compact. Indeed, let
\[
\zeta=
\min_{\substack{\sigma\in \mathscr{P}_r\\1\leq l\leq r-1}}
|s_l^\sigma|^p.
\]
It follows from the CNC that $\zeta>0$. Moreover, for any $v,w\in T$, we have
\[
A_D(v,w)=\chi_{\{w\}}(\varphi_i(v))-\chi_{\{w\}}(\varphi_j(v)).
\]
Whenever $A_D(v,w)=0$, the inequality \[|A_Q(v,w)|^p\geq \zeta|A_D(v,w)|^p\] trivially holds.
If instead $A_D(v,w)\ne0$, then exactly one of $\varphi_i(v)$ and $\varphi_j(v)$ equals $w$.
Put $I=\{k\in\{1,\ldots,r\}:\varphi_k(v)=w\}$. It follows that $1\leq|I|\leq r-1$.
We may choose $\sigma=(\sigma_1,\ldots,\sigma_r)\in \mathscr{P}_r$ such that the first $|I|$ entries
$\sigma_1,\ldots,\sigma_{|I|}$ are precisely the elements of $I$.
Thus we obtain that
\begin{align*}
A_Q(v,w)&=\sum_{k=1}^{r}a_k\chi_{\{w\}}(\varphi_k(v))\\
&=\sum_{j=1}^{|I|}a_{\sigma_j}=s_{|I|}^{\sigma}.
\end{align*}
Notice that $|A_D(v,w)|=1$. It follows that
\[
|A_Q(v,w)|^p\geq \zeta|A_D(v,w)|^p
\]
for all $v,w\in T$. Without loss of generality, we may assume $v\in D_n$ and $w\in D_m$. Summing over all $v\in D_n$, we arrive at
\begin{align*}
N_Q(n,w)
&=\sum_{|v|=n}|A_Q(v,w)|^p\\
&\geq
\zeta\sum_{|v|=n}|A_D(v,w)|^p\\
&=\zeta N_D(n,w),
\end{align*}
and hence
\[
\max_{|w|=m}N_D(n,w)
\leq
\frac1{\zeta}
\max_{|w|=m}N_Q(n,w).
\]
Consequently, by Theorem \ref{thm}, we have
\begin{align*}
0&\leq
\sup_{n\in\mathbb N_0}\frac1{c_n}\sum_{m>R}c_m\max_{|w|=m}N_D(n,w)\\
&\leq\frac1{\zeta}\sup_{n\in\mathbb N_0}\frac1{c_n}\sum_{m>R}c_m\max_{|w|=m}N_Q(n,w)\to0
\end{align*}
as $R\to\infty$, where the convergence follows from the compactness of $Q$. Using Theorem \ref{thm} again, we conclude that $D=C_{\varphi_i}-C_{\varphi_j}$ is compact on $\mathbb T_p$. The proof is complete.
\end{proof}

\begin{example}
Theorem \ref{thm3.5} allows us to readily characterize the compactness of the CNC admissible linear combinations of composition operators on $\mathbb T_p$.
\begin{itemize}
\item The operator $9C_{\varphi_1}-5C_{\varphi_2}+C_{\varphi_3}-3C_{\varphi_4}$ is compact on $\mathbb T_p$ if and only if $C_{\varphi_i}$, $i=1,...4$, are all compact on $\mathbb T_p$.
\item The operator $2C_{\varphi_1}-C_{\varphi_2}-C_{\varphi_3}$ is compact on $\mathbb T_p$ if and only if $C_{\varphi_i}-C_{\varphi_j}$ is compact for each $i,j=1,2,3$.
\end{itemize}
\end{example}

Let $\varphi, \psi$ be holomorphic self-maps of the unit disk $\mathbb{D}$ in the complex plane $\mathbb{C}$. For $0<p<\infty$, let $H^{p}(\mathbb{D})$ stand for the Hardy space on $\mathbb{D}$.
In \cite{cckp}, Choe et al. formulated the following conjecture.

\begin{conjecture}\label{conj}
\cite[Conjecture 4.3]{cckp}
The compactness of $C_\varphi-C_\psi$ on $H^{p}(\mathbb{D})$ is independent of the parameter $p\in(0, \infty)$.
\end{conjecture}

Applying Theorem \ref{thm}, we deduce the following result.

\begin{theorem}\label{thm3.8}
Let $T$ be a $(q+1)$-homogeneous tree with $q\ge2$, and let $1\le p<\infty$. For each $i=1, \ldots, r$, let $\varphi_i: T\to T$ be a self-map such that the composition operator $C_{\varphi_i}$ is bounded on $\mathbb T_p$. Then the compactness of the operator $Q=\sum_{i=1}^{r} a_iC_{\varphi_i}$ (where $0\neq a_i\in\mathbb C$) on $\mathbb T_{p}$ is independent of the parameter $p\in [1, \infty)$. In other words, for all $1\leq p_1\neq p_2<\infty$,
\[
Q\text{ is compact on }\mathbb T_{p_1}.
\quad\Longleftrightarrow\quad
Q\text{ is compact on }\mathbb T_{p_2}.
\]
\end{theorem}

\begin{proof}
Without loss of generality, suppose $p_1<p_2$. Set $\Lambda_r=\{1,\ldots,r\}$. For any $v,w\in T$, let
\[
I(v,w)=\{i\in\Lambda_r:\varphi_i(v)=w\}.
\]
By the definition of $A_Q$, we have
\[
A_Q(v,w)=\sum_{i\in I(v,w)}a_i.
\]
Moreover, we define
\[
C_1=\min\left\{\left|\sum_{i\in I}a_i\right|:
I\subset\Lambda_r,\ \sum_{i\in I}a_i\neq0
\right\}
\]
and
\[
C_2=\max\left\{\left|\sum_{i\in I}a_i\right|:
I\subset\Lambda_r
\right\}.
\]
It follows that
\[
0<C_1\leq C_2<\infty.
\]
Whenever $A_Q(v,w)\neq0$, we have
\[
C_1\leq |A_Q(v,w)|\leq C_2,
\]
and hence
\[
C_1^{p_2-p_1}\leq |A_Q(v,w)|^{p_2-p_1}\leq C_2^{p_2-p_1}.
\]
This yields
\[
C_1^{p_2-p_1}|A_Q(v,w)|^{p_1}
\leq |A_Q(v,w)|^{p_2}
\leq C_2^{p_2-p_1}|A_Q(v,w)|^{p_1}.
\]
If $A_Q(v,w)=0$, then the above inequality holds trivially.		
Summing the preceding inequality over all $v\in D_n$, we obtain
\[
C_1^{p_2-p_1}N_{Q,p_1}(n,w)
\leq N_{Q,p_2}(n,w)
\leq C_2^{p_2-p_1}N_{Q,p_1}(n,w),
\]
from which it follows that
\[
C_1^{p_2-p_1}\max_{|w|=l}N_{Q,p_1}(n,w)
\leq \max_{|w|=l}N_{Q,p_2}(n,w)
\leq C_2^{p_2-p_1}\max_{|w|=l}N_{Q,p_1}(n,w),
\]
where $N_{Q,p_j}(n,w)=\sum_{|v|=n}|A_Q(v,w)|^{p_j}$ ($j=1,2$).
Consequently, for each $R\in\mathbb N_0$, we have
\[
\begin{aligned}
C_1^{p_2-p_1}
&\sup_{n\in\mathbb N_0}\frac1{c_n}
\sum_{l>R}c_l\max_{|w|=l}N_{Q,p_1}(n,w)\\
&\leq
\sup_{n\in\mathbb N_0}\frac1{c_n}
\sum_{l>R}c_l\max_{|w|=l}N_{Q,p_2}(n,w)\\
&\leq
C_2^{p_2-p_1}
\sup_{n\in\mathbb N_0}\frac1{c_n}
\sum_{l>R}c_l\max_{|w|=l}N_{Q,p_1}(n,w).
\end{aligned}
\]
Since $0<C_1\leq C_2<\infty$, we deduce from the above inequality that
\[
\begin{aligned}
&\lim_{R\to\infty}
\sup_{n\in\mathbb N_0}\frac1{c_n}
\sum_{l>R}c_l\max_{|w|=l}N_{Q,p_1}(n,w)=0\\
&\qquad\Longleftrightarrow
\lim_{R\to\infty}
\sup_{n\in\mathbb N_0}\frac1{c_n}
\sum_{l>R}c_l\max_{|w|=l}N_{Q,p_2}(n,w)=0.
\end{aligned}
\]
By Theorem \ref{thm}, we further obtain that
\[
Q\text{ is compact on }\mathbb T_{p_1}.
\quad\Longleftrightarrow\quad
Q\text{ is compact on }\mathbb T_{p_2}.
\]
\end{proof}

In particular, by Theorem \ref{thm3.8}, we immediately obtain the following corollary, which gives an affirmative answer to Conjecture \ref{conj} for discrete structures when $1\leq p<\infty$.

\begin{corollary}
Let $T$ be a $(q+1)$-homogeneous rooted tree with $q\geq 2$ and $1\leq p<\infty$. Let $\varphi, \psi$ be self-maps of $T$ such that $C_{\varphi}$ and $C_{\psi}$ are bounded. Then the compactness of the operator  $C_{\varphi}-C_{\psi}$ on $\mathbb{T}_p$ does not depend on the parameter $p\in [1, \infty)$.
\end{corollary}

\section{Examples}

In this section we construct two examples to demonstrate that Moorhouse cancellation property for compact (double) differences of composition operators is valid on the space $\mathbb{T}_p$.

The first example shows that a pair of noncompact composition operators can give rise to a compact difference on $\mathbb{T}_p$.

\begin{example}
There exist self-maps $\varphi,\psi: T\to T$ satisfying the following properties:
\begin{enumerate}[(i)]
\item $C_\varphi$ and $C_\psi$ are bounded on $\mathbb T_p$;
\item Neither $C_\varphi$ nor $C_\psi$ is compact on $\mathbb T_p$;
\item $C_\varphi-C_\psi$ is compact on $\mathbb T_p$.
\end{enumerate}
\end{example}

\begin{proof}		
Choose a sufficiently large integer $N_0\in\mathbb{N}_0$ such that $c_{\lfloor n/2\rfloor}\ge 2$ holds for all $n\ge N_0$. For every $n\ge N_0$, select a vertex $v_n\in D_n$ together with two distinct vertices $x_n, y_n\in D_{\lfloor n/2\rfloor}$. Define the self-maps $\varphi$ and $\psi$ by
\[
\varphi(v)=
\begin{cases}
x_n, & v=v_n,\ n\ge N_0,\\
v, & v\notin\{v_n:n\ge N_0\},
\end{cases}
\qquad
\psi(v)=
\begin{cases}
y_n, & v=v_n,\ n\ge N_0,\\
v, & v\notin\{v_n:n\ge N_0\}.
\end{cases}
\]

We first prove that $C_\varphi$ and $C_\psi$ are bounded on $\mathbb T_p$.
Indeed, for any $f\in\mathbb T_p$ and $n\geq N_0$, we have
\begin{align*}
M_p^p(n,C_\varphi f)
&=\frac1{c_n}\sum_{|v|=n}|f(\varphi(v))|^p\\
&=\frac1{c_n}
\sum_{\substack{|v|=n\\v\neq v_n}}|f(v)|^p
+\frac1{c_n}|f(x_n)|^p.
\end{align*}
Since $x_n\in D_{\lfloor n/2\rfloor}$, it follows that
\[
|f(x_n)|^p\leq c_{\lfloor n/2\rfloor}\|f\|^p .
\]
Therefore,
\begin{align*}
M_p^p(n,C_\varphi f)
&\leq
\|f\|^p+
\frac{c_{\lfloor n/2\rfloor}}{c_n}\|f\|^p\\
&\leq 2\|f\|^p .
\end{align*}
For any $0\leq n<N_0$, it is clear that
$M_p(n,C_\varphi f)=M_p(n,f)$. Hence
\[
\|C_\varphi f\|\leq 2^{1/p}\|f\|,
\]
which implies that $C_\varphi$ is bounded on $\mathbb T_p$.
Similarly, $C_\psi$ is also bounded on $\mathbb T_p$.

Next, we prove that neither $C_\varphi$ nor $C_\psi$ is compact on $\mathbb T_p$. Indeed, for every
$n\geq N_0$ and $w\in D_n\setminus\{v_n\}$, we have
\[
N_\varphi(n,w)\geq1.
\]
It follows that
\[
\max_{|w|=n}N_\varphi(n,w)\geq1 .
\]
For any $R\in\mathbb N_0$, take $n>\max\{R,N_0\}$. Then
\begin{align*}
\sup_{s\in\mathbb N_0}
\frac1{c_s}\sum_{m>R}c_m
\max_{|w|=m}N_\varphi(s,w)
&\geq
\frac1{c_n}c_n
\max_{|w|=n}N_\varphi(n,w)\\
&\geq1.
\end{align*}
Consequently,
\[
\sup_{s\in\mathbb N_0}
\frac1{c_s}\sum_{m>R}c_m
N_{m,s}\nrightarrow0
\]
as $R\to\infty$. By Corollary \ref{cor3.4}, $C_\varphi$ is not compact on $\mathbb T_p$.
Similarly, $C_\psi$ is not compact on $\mathbb T_p$.

Finally, we prove that $P:=C_\varphi-C_\psi$ is compact. Indeed, for any $v, w\in T$, we have
\[
A_P(v,w)=\chi_{\{w\}}(\varphi(v))-\chi_{\{w\}}(\psi(v)).
\]
If $v\neq v_n$ for all $n\geq N_0$, then
$\varphi(v)=\psi(v)=v$, so that
\[
A_P(v,w)=0.
\]
If $v=v_n$ for some $n\geq N_0$, then
\[
\varphi(v_n)=x_n,\qquad \psi(v_n)=y_n,
\]
and hence
\[
A_P(v_n,w)=\chi_{\{w\}}(x_n)-\chi_{\{w\}}(y_n).
\]
It follows that
\[
N_P(n,w)=
\begin{cases}
1, & w=x_n\ \text{or}\ w=y_n,\\
0, & \text{otherwise}.
\end{cases}
\]
Since $x_n, y_n\in D_{\lfloor n/2\rfloor}$, we obtain that
\[
\max_{|w|=m}N_P(n,w)=
\begin{cases}
1, & m=\lfloor n/2\rfloor,\\
0, & \text{otherwise}.
\end{cases}
\]
Therefore,
\begin{align*}
\sup_{n\in\mathbb N_0}
\frac1{c_n}\sum_{m>R}c_m
\max_{|w|=m}N_P(n,w)=
\sup_{\substack{n\geq N_0\\\lfloor n/2\rfloor>R}}
\frac{c_{\lfloor n/2\rfloor}}{c_n}.
\end{align*}
Set $k=\lfloor n/2\rfloor$. If $n=2k$, then
\[
\frac{c_{\lfloor n/2\rfloor}}{c_n}
=\frac{c_k}{c_{2k}}
=q^{-k}.
\]
If $n=2k+1$, then
\[
\frac{c_{\lfloor n/2\rfloor}}{c_n}
=\frac{c_k}{c_{2k+1}}
=q^{-(k+1)} .
\]
Recalling that $k>R$, it follows in both cases that
\[
\frac{c_{\lfloor n/2\rfloor}}{c_n}
\leq q^{-(R+1)} .
\]
Hence
\[
\sup_{n\in\mathbb N_0}
\frac1{c_n}\sum_{m>R}c_m
\max_{|w|=m}N_P(n,w)
\leq q^{-(R+1)}\rightarrow0
\]
as $R\to\infty$.
Using Theorem \ref{thm}, we see that $P=C_\varphi-C_\psi$ is compact on $\mathbb T_p$.
\end{proof}

The example below illustrates that double difference cancellation can occur for linear combinations of four composition operators acting on $\mathbb T_p$.

\begin{example}
There exist self-maps $\varphi_i: T\to T$, $i=1,...,4$ satisfying the following properties:
\begin{enumerate}[(i)]
\item $C_{\varphi_i}$, $i=1,...,4$, are bounded but not compact on $\mathbb T_p$;
\item None of $C_{\varphi_i}-C_{\varphi_j}$ ($i\neq j$) is compact on $\mathbb T_p$;
\item $(C_{\varphi_1}-C_{\varphi_2})-(C_{\varphi_3}-C_{\varphi_4})$ is compact on $\mathbb T_p$.
\end{enumerate}
\end{example}

\begin{proof}		
Set
\[
E=\bigcup_{l\geq0}D_{3l+1},
\qquad
F=\bigcup_{l\geq0}D_{3l+2}.
\]
For each $k\in\mathbb N$, pick a vertex $u_k\in D_{3k}$ and two distinct vertices $x_k, y_k\in D_k$. Let
$\varphi_i(o)=o$, $i=1,2,3,4$. For $v\neq o$, define
\[
\varphi_1(v)=
\begin{cases}
v, & v\in E\cup F,\\
v^{-}, & v\in D_{3k}\setminus\{u_k\},\\
y_k, & v=u_k,
\end{cases}
\qquad
\varphi_2(v)=
\begin{cases}
v, & v\in E,\\
v^{-}, & v\in F,\\
v^{-}, & v\in D_{3k}\setminus\{u_k\},\\
y_k, & v=u_k,
\end{cases}
\]
and
\[
\varphi_3(v)=
\begin{cases}
v^{-}, & v\in E,\\
v, & v\in F,\\
v^{-}, & v\in D_{3k}\setminus\{u_k\},\\
y_k, & v=u_k,
\end{cases}
\qquad
\varphi_4(v)=
\begin{cases}
v^{-}, & v\in E\cup F,\\
v^{-}, & v\in D_{3k}\setminus\{u_k\},\\
x_k, & v=u_k.
\end{cases}
\]
	
We first show that the operators $C_{\varphi_j}$, $j=1,2,3,4$, are bounded on $\mathbb T_p$.
By \cite[Theorem 4]{mp}, it suffices to verify
\[
\sup_{n\in\mathbb N_0}
\frac1{c_n}\sum_{m=0}^{\infty}c_mN_{m,n}<\infty
\]
for all $\varphi_j$, $j=1,2,3,4$.
Indeed, for $n=3k+1$ with $k\geq0$, it follows from the definitions of $\varphi_1$ and $\varphi_2$ that $N_{n,n}=1$ and $N_{m,n}=0$ whenever $m\neq n$. Therefore
\[
\frac1{c_n}\sum_{m=0}^{\infty}c_mN_{m,n}=1.
\]
By the definitions of $\varphi_3$ and $\varphi_4$, we have $\varphi_3(v)=\varphi_4(v)=v^{-}$ for all $v\in D_n$. If $n>1$, then $N_{n-1,n}=q$ and $N_{m,n}=0$ for all $m\neq n-1$. It follows that
\[
\frac1{c_n}\sum_{m=0}^{\infty}c_mN_{m,n}
=\frac{qc_{n-1}}{c_n}=1.
\]
If $n=1$, then $N_{0,1}=c_1$, and the corresponding quantity equals to $1$.
		
For $n=3k+2$ with $k\geq0$, $\varphi_1$ and $\varphi_3$ act as the identity on $D_n$, while
$\varphi_2$ and $\varphi_4$ coincide with the parent map on $D_n$.
The same calculation as above yields
\[
\frac1{c_n}\sum_{m=0}^{\infty}c_mN_{m,n}=1
\]
for each of the four maps.

For $n=3k$ with $k>0$, by the definitions of $\varphi_i$, $i=1,...,4$, we have $N_{3k-1,3k}=q$ and $N_{k,3k}=1$.
It follows that
\[
\begin{aligned}
\frac1{c_{3k}}\sum_{m=0}^{\infty}c_mN_{m,3k}
&=\frac{qc_{3k-1}+c_k}{c_{3k}}\\
&=1+\frac{c_k}{c_{3k}}
=1+q^{-2k}
\leq1+q^{-2}.
\end{aligned}
\]

Hence we conclude that
\[
\sup_{n\in\mathbb N_0}
\frac1{c_n}\sum_{m=0}^{\infty}c_mN_{m,n}
\leq1+q^{-2}<\infty,
\]
which yields that $C_{\varphi_i}$, $j=1,...,4$, are bounded on $\mathbb T_p$.
		
Now we show that the operators $C_{\varphi_j}$, $j=1,2,3,4$, are not compact on $\mathbb T_p$ by applying Corollary \ref{cor3.4}.
Indeed, for each $k\geq0$ and $w\in D_{3k+1}$, we have $\varphi_i(w)=w$ ($i=1,2$). Then
$N_{\varphi_i}(3k+1,w)=1$, and hence
\[
\max_{|w|=3k+1}N_{\varphi_i}(3k+1,w)=1
\]
for $i=1,2$. Given any $R\in\mathbb N_0$, there exists an integer $k$ with $3k+1>R$. Thus
\[
\begin{aligned}
\sup_{n\in\mathbb N_0}
\frac1{c_n}\sum_{m>R}c_m
\max_{|w|=m}N_{\varphi_1}(n,w)
&\geq
\frac1{c_{3k+1}}c_{3k+1}
\max_{|w|=3k+1}N_{\varphi_1}(3k+1,w)\\
&=1.
\end{aligned}
\]
Using Corollary \ref{cor3.4}, we conclude that $C_{\varphi_i}$ ($i=1,2$) are not compact.
		
Similarly, for each $k\geq0$ and $w\in D_{3k+2}$, we have $\varphi_3(w)=w$. Hence the same reasoning as above yields that $C_{\varphi_3}$ is not compact.
		
Moreover, for each $k>0$ and $v\in D_{3k+1}$, we have $\varphi_4(v)=v^{-}$. Then
\[
\max_{|w|=3k}N_{\varphi_4}(3k+1,w)=q.
\]
For every $R\in\mathbb N_0$, we may choose $k$ such that $3k>R$. It follows that
\[
\begin{aligned}
\sup_{n\in\mathbb N_0}
\frac1{c_n}\sum_{m>R}c_m
\max_{|w|=m}N_{\varphi_4}(n,w)
&\geq
\frac{c_{3k}}{c_{3k+1}}
\max_{|w|=3k}N_{\varphi_4}(3k+1,w)\\
&=\frac{qc_{3k}}{c_{3k+1}}
=1.
\end{aligned}
\]
By Corollary \ref{cor3.4}, $C_{\varphi_4}$ fail to be compact. Therefore, property (i) holds.
		
Next, we establish property (ii). Set
\[
D_{1,2}=C_{\varphi_1}-C_{\varphi_2}.
\]
For every $k\geq 0$ and $w\in D_{3k+2}$, we have $\varphi_1(w)=w$ and $\varphi_2(w)=w^{-}$. This gives $A_{D_{1,2}}(w,w)=1$. On the other hand, for any $v\in D_{3k+2}$ with $v\neq w$, we have $\varphi_1(v)\neq w$ and $\varphi_2(v)\neq w$. It follows that
\[
N_{D_{1,2}}(3k+2,w)=1
\]
for all $w\in D_{3k+2}$. As a result,
\[
\max_{|w|=3k+2}N_{D_{1,2}}(3k+2,w)=1.
\]
For each $R\in\mathbb N_0$, there exists an integer $k$ such that $3k+2>R$. Then
\[
\begin{aligned}
\sup_{n\in\mathbb N_0}
\frac1{c_n}\sum_{m>R}c_m
\max_{|w|=m}N_{D_{1,2}}(n,w)
&\geq
\frac1{c_{3k+2}}c_{3k+2}
\max_{|w|=3k+2}N_{D_{1,2}}(3k+2,w)\\
&=1.
\end{aligned}
\]
By Theorem \ref{thm}, we deduce that $D_{1,2}=C_{\varphi_1}-C_{\varphi_2}$ is not compact.

Similarly, one can verify that the operators $C_{\varphi_3}-C_{\varphi_4}$, $C_{\varphi_1}-C_{\varphi_3}$ and $C_{\varphi_2}-C_{\varphi_4}$ all fail to be compact.
		
Finally, we show that the operator $Q:=(C_{\varphi_1}-C_{\varphi_2})-(C_{\varphi_3}-C_{\varphi_4})$ is compact. Indeed, for any $v,w\in T$, we have
\[
\begin{aligned}
A_Q(v,w)
&=\chi_{\{w\}}(\varphi_1(v))
-\chi_{\{w\}}(\varphi_2(v))\\
&\quad-\chi_{\{w\}}(\varphi_3(v))
+\chi_{\{w\}}(\varphi_4(v)).
\end{aligned}
\]
If $v\in D_{3k+1}$ with $k\geq0$, then
\[
\varphi_1(v)=\varphi_2(v)=v, \qquad \varphi_3(v)=\varphi_4(v)=v^{-},
\]
which implies $A_Q(v,w)=0$. Similarly, if $v\in D_{3k+2}$ with $k\geq0$, then
\[
\varphi_1(v)=\varphi_3(v)=v, \qquad \varphi_2(v)=\varphi_4(v)=v^{-},
\]
so that $A_Q(v,w)=0$. If $v\in D_{3k}\setminus\{u_k\}$ with $k>0$, then
\[
\varphi_1(v)=\varphi_2(v)=\varphi_3(v)=\varphi_4(v)=v^{-},
\]
and hence $A_Q(v,w)=0$.

Furthermore, if $v=u_k$, then we have
\[
\varphi_1(u_k)=\varphi_2(u_k)=\varphi_3(u_k)=y_k,
\qquad
\varphi_4(u_k)=x_k.
\]
Hence
\[
A_Q(u_k,w)
=-\chi_{\{w\}}(y_k)+\chi_{\{w\}}(x_k).
\]
Since $x_k\neq y_k$, it follows that
\[
N_Q(3k,w)=
\begin{cases}
1, & w=x_k\ \text{or}\ w=y_k,\\
0, & \text{otherwise}.
\end{cases}
\]
Notice that $x_k, y_k\in D_k$. Thus,
\[
\max_{|w|=m}N_Q(3k,w)=
\begin{cases}
1, & m=k,\\
0, & m\neq k.
\end{cases}
\]		
Therefore, for any $R\in\mathbb N_0$, we obtain that
\[
\begin{aligned}
&\sup_{n\in\mathbb N_0}
\frac1{c_n}\sum_{m>R}c_m
\max_{|w|=m}N_Q(n,w)\\
&\qquad=\sup_{k>R}\frac{c_k}{c_{3k}}
=\sup_{k>R}q^{-2k}
=q^{-2(R+1)}.
\end{aligned}
\]
Hence
\[
\sup_{n\in\mathbb N_0}
\frac1{c_n}\sum_{m>R}c_m
\max_{|w|=m}N_Q(n,w)
\to 0
\]
as $R\to\infty$. By Theorem \ref{thm}, we conclude that
$(C_{\varphi_1}-C_{\varphi_2})-(C_{\varphi_3}-C_{\varphi_4})$
is compact on $\mathbb T_p$.
\end{proof}

\begin{remark}
The Moorhouse cancellation phenomenon for compact double differences generated by three composition operators cannot occur on $\mathbb T_p$. Namely, if $(C_{\varphi_1}-C_{\varphi_2})-(C_{\varphi_3}-C_{\varphi_1})$ is compact on $\mathbb T_p$, then by Theorem \ref{thm3.5} (ii), both $C_{\varphi_1}-C_{\varphi_2}$ and $C_{\varphi_3}-C_{\varphi_1}$ are compact.
\end{remark}

{\bf Acknowledgement} This work is supported by the Natural Science Basic Research Program of Shaanxi Province (Program No. 2025JC-YBMS-024).

{\bf Author Contributions} All authors share credit for this work equally.

{\bf Data Availability} No datasets were generated or analysed during the current study.

{\bf Declarations}

{\bf Conflict of interest} The authors declare no competing interests.

\bibliographystyle{amsplain}

\end{document}